\documentclass[12pt]{article}

\usepackage{erics_preprints}

\input xy
\xyoption{all}
\xyoption{poly}
\xyoption{matrix}

\begin{document}
\title{On a Generalization of the Notion of Semidirect Product of Groups}
\author{Eric Ram\'on Antokoletz}
\date{}
\maketitle
\begin{abstract}
We introduce an {\it external} version of the 
{\it internal $r$-fold semidirect product of groups
(SDP)} of
\cite{CarrascoCegarra}. Just as for the classical external SDP, 
certain algebraic data are required to guarantee associativity
of the construction. We give an algorithmic procedure for computing
axioms characterizing these data. 
Additionally, we give criteria for determining when
a family of homomorphisms from the factors of an SDP into a monoid or
group assemble into a homomorphism on the entire SDP. These tools will
be used elsewhere to give explicit algebraic axioms for {\it
hypercrossed complexes}, which are algebraic models for classical
homotopy types introduced in \cite{CarrascoCegarra}.
\end{abstract}
\tableofcontents
	
\setcounter{secnumdepth}{1}
\setcounter{section}{0}
							
	\section{Introduction}							
	\label{Intro_SDPs}


In this paper we lay foundations for {\it $r$-fold semidirect products
of groups}, which constitute a particular generalization of the usual
notion of semidirect product of two groups. This generalization was
introduced by P. Carrasco and A. M. Cegarra in \cite{CarrascoCegarra}
and shown there to play a fundamental role in the structure of
simplicial groups. More specifically, it enabled them to state a
nonabelian Dold-Kan theorem to the effect that a simplicial group $G$ is
equivalent to a certain nonabelian differential-algebraic structure,
dubbed a {\it hypercrossed complex} in \cite{CarrascoCegarra},
consisting of the Moore complex of $G$ together with a large family of
pairings among the terms of the Moore complex. (The reader will recall
that the Moore complex $M(G)$ of a simplicial group $G$, named after its
discoverer John C. Moore, is a chain complex of nonabelian groups whose
homology groups are the homotopy groups of $G$.  See \cite{JCMoore} or
\cite{May}). 

This is of interest because, according to an old result of Daniel Kan,
the homotopy theory of simplicial groups is equivalent to the homotopy
theory of pointed connected topological spaces (simplicial groups play the role
of loop spaces---see \cite{KanCombinHomotopy} and \cite{KanLine} for
Kan's original results and also \cite{Quillen} and \cite{Hovey} for
their modern reformulation in terms of model categories). Therefore the
aforementioned result of Carrasco and Cegarra can be said to place the
homotopy theory of spaces in a homological-algebraic context---namely
the category of hypercrossed complexes---just as earlier the Moore
complex by itself had done so for the homotopy groups of spaces. 

In order to delineate the role in all of this of Carrasco and Cegarra's
$r$-fold semidirect products, we proceed to describe their result. The
classical Dold-Kan theorem (\cite{May}, \cite{GoerssJardine}) shows that
the Moore complex functor $M$, restricted to the category ${\bf SAb}$ of
simplicial abelian groups, is already an equivalence of categories
between ${\bf SAb}$ and the category ${\bf Ch}_+({\mathbb Z})$ of
nonnegative chain complexes of abelian groups. Inherent in this
equivalence is the fact that the chain complex $M(A)$ by itself uniquely
determines (up to isomorphism) the simplicial abelian group $A$ from
which it originated. Briefly, each term $A_n$ of $A$ decomposes as a
direct sum, indexed by simplicial degeneracy operators, of the terms
$M_k(A)$ of the chain complex $M(A)$, and one can reconstruct the
simplicial operators on $A$ by using the simplicial identities as a
formal bookkeeping system, with occasional input from the boundary
operator of $M(A)$.

The general case of an arbitrary (nonabelian) simplicial group $G$ is
not so simple, as the Moore complex by itself is in general not enough
to determine $G$. Closely related is the fact that, although the $n$th
homotopy groups of simplicial groups are well-known to be abelian for
$n\ge 1$ (see \cite{May} or \cite{Lamotke}), one may easily construct
chain complexes of nonabelian groups whose higher homology groups are
not abelian. Thus one can see that there are extra constraints on Moore
complexes distinguishing them from among arbitrary chain complexes of
nonabelian groups. This situation leads to the following question: for a
given arbitrary simplicial group $G$, what extra information must be
retained along with the Moore complex $M(G)$ in order to completely
reconstruct $G$ from $M(G)$? The complex $M(G)$ together with this extra
information is precisely the definition of {\it hypercrossed complex} a
la Carrasco and Cegarra.

The first step towards answering this question is to identify the
correct analog of direct sum appearing in the abelian Dold-Kan theorem.
Direct sums appear there because surjective morphisms in the simplicial
category $\Delta$ possess right inverses, so that every degeneracy
operator on a simplicial abelian group gives a direct sum decomposition
by embedding its source as a direct summand of its target. In the
nonabelian case, the most that can be said in general is that every
degeneracy operator $s_{i_k}\dots s_{i_1}:G_{n-k}\ra G_n$ on a
simplicial group $G$ embeds its source as a {\it retract} of its target.
Accordingly (\cite{CarrascoCegarra},\cite{preprint2}), in place of
direct sums, the nonabelian Dold-Kan theorem of Carrasco and Cegarra
states that, for each $n\ge 0$, a certain family of intersections of the
kernels of these retractions on $G_n$ can be arranged into a filtration
of $G_n$ by normal subgroups which is completely split in the sense that
each subquotient is canonically isomorphic to a copy of a term
$M_{n-k}(G)$ of the Moore complex embedded in $G_n$ via the degeneracy
$s_{i_k}\dots s_{i_1}:G_{n-k}\ra G_n$. This type of decomposition is
then what Carrasco and Cegarra identify as the correct generalization in
this context of the notion of {\it internal semidirect product} (see section
\ref{InternalSDPs} for alternative formulations of this definition,
including the original one introduced in \cite{CarrascoCegarra}).

The next step, which is the goal of the present paper, is to provide an
axiomatic description of {\it external $r$-fold semidirect products},
that is, an axiomatization of the data required to reconstruct an
$r$-fold semidirect product from its factors. These axioms can then be
applied to the Dold-Kan decomposition of \cite{CarrascoCegarra} to
enable a characterization of the pairings appearing in hypercrossed
complexes (this has been done by the author and will appear in a
forthcoming article).


These generalized semidirect products can moreover be expected to play a similar
role in a wider class of cases of interest. In \cite{GenReedy}, Ieke Moerdijk and
Clemens Berger introduce the notion of {\it Eilenberg-Zilber category}
$R$ which is a special case of their {\it generalized Reedy categories}
and includes such examples as the cubical category $\Box$ as well as the
cyclic category $\Lambda$ and more generally the total category of any
crossed simplicial group (see \cite{GenReedy} for definitions and even more general examples).
These categories all share with $\Delta$ the property mentioned above,
that epimorphisms admit right inverses, and so group-valued functors
coming from these categories can also be expected to admit nonabelian
Dold-Kan theorems in the spirit of \cite{CarrascoCegarra} and \cite{preprint2}.

	\section{Internal $r$-fold Semidirect Products (SDPs)}			

	\label{InternalSDPs}

In this section, we give several equivalent formulations of the
definition of semidirect product, including the original one introduced
in \cite{CarrascoCegarra}.

\bigskip

Assume given a group $G$ and subgroups $H_1,\ldots,H_r$.

\defn{ \label{iSDP}
	The group $G$ is said to be an {\it internal $r$-semidirect product
	\emph{(briefly an {\it SDP})}~of the subgroups $H_i$} if the
	following two conditions hold.
	\begin{enumerate}
		\item{The set $H_1H_2\ldots H_i$ is a normal subgroup of $G$ for all $i$.}
		\item{Every $g\in G$ can be factored uniquely as a product 
			$$g = h_1h_2\ldots h_r$$ with $h_i\in H_i$ for all $i$.}
	\end{enumerate}
	We follow \cite{CarrascoCegarra} in using the notation
		$$G = H_1\rtimes\ldots\rtimes H_r$$
	if the above conditions hold.
}

\eg{	\label{iDP}
	If condition 1 is replaced by the stronger condition 
		$$H_i \lhd G \mforall i$$
	then the definition reduces to that of internal {\it direct} product, which is
	thus a special case of an internal SDP.
}

\rmk{	The order in which the subgroups $H_i$ appear is an essential part of
	the definition of SDP. It can happen that $G$ is an SDP of the $H_i$
	when they arranged in certain orders, but not in others. At one extreme,
	there may be a unique such order. At the other extreme, it is
	immediate that $G$ is an internal direct product of the
	$H_i$ if and only if $G$ is an SDP of the $H_i$ arranged in each
	possible order.
}

The following proposition contains the original definition given in
\cite{CarrascoCegarra} and demonstrates that Definition \ref{iSDP} above is
equivalent to it.

\begin{prp}
	The group $G$ is an internal $r$-SDP of the subgroups $H_i$ if and
	only if the following three conditions hold.
	\emph{\begin{enumerate}
		\item{The set $H_1H_2\ldots H_i$ is a normal subgroup of $G$ for all $i$.}
		\item[$\mathrm{2a}$.]{$G = H_1H_2\ldots H_r$}
		\item[$\mathrm{2b}$.]{$H_1H_2\ldots H_i\cap H_{i+1} = \{1\}$ for $1\leq i \leq r-1$}
	\end{enumerate}}
\end{prp}
\prf{Proof.}{
	It will be shown that, in the presence of condition 1, condition 2
	is equivalent to the conjunction of 2a and 2b.
	For this it suffices to assume something weaker than 1, 
	namely that the sets $H_1H_2\ldots H_i$ are merely
	subgroups of $G$ for all $i$.

	Condition 2a is equivalent to existence of the factorizations of
	condition 2, and it remains to show that 2b is equivalent to uniqueness of the
	factorizations. Uniqueness implies 2b since any nontrivial element of
	$H_1H_2\ldots H_i\cap H_{i+1}$ would evidently possess two distinct
	factorizations. Conversely, assume 2b and let $g\in G$ have the
	following two factorizations.
		$$g = h_1\ldots h_r = h_1^\p \ldots h_r^\p$$
	It follows that
		$$(h_1^\p \ldots h_{r-1}^\p)^{-1}(h_1\ldots h_{r-1}) = h_r^\p h_r^{-1}$$
	and since by hypothesis $H_1\ldots H_{r-1}$ is a subgroup, the left hand
	side belongs to it. Condition 2b now implies
		$$h_r = h_r^\p$$
		$$h_1\ldots h_{r-1} = h_1^\p \ldots h_{r-1}^\p$$
	and, repeating the argument inductively, one deduces that the two
	factorizations of $g$ coincide. 
}

For the next alternate formulation of the definition, the following common notations are used.  
For elements $h,k\in G$, let 
	$${}^kh := khk^{-1}$$
	$$[k,h] := khk^{-1}h^{-1}$$
denote respectively the {\it conjugate of $h$ by $k$} and the {\it commutator of $k$ and $h$}.
For subgroups $H,K$ of a group $G$, let
	$${}^KH := \setof{ {}^kh }{h\in H, k\in K}$$
denote the subset of $G$ consisting of all conjugates of elements of $H$ by elements of $K$.
The familiar notation $[K,H]$ is also used
for the subgroup of $G$ generated by all commutators $[k,h]$
with $h\in H, k\in K$.

\bigskip

The following proposition enables condition 1 of the definition to be checked ``term by term''.

\begin{prp}  \label{checktermbyterm}
	Provided that $G$ is generated by the subgroups $H_i$,
	condition $1$ is equivalent to each of the following.
	\emph{\begin{enumerate}
		\item[$1^\p$.]{${}^{H_j}H_i \subseteq H_1\ldots H_i$ for any $i,j$ with $1\le i < j \le r$.}
		\item[$1^\pp$.]{$[H_j,H_i] \subseteq H_1\ldots H_i$ for any $i,j$ with $1\le i < j \le r$.}
	\end{enumerate}}
\end{prp}
\prf{Proof.}{
	The equivalence of $1^\p$ and $1^\pp$ follows immediately from
	the identity
		$$[k,h] = {}^khh^{-1}$$
	so it suffices to prove that condition 1 is equivalent to $1^\p$.
	It is also immediate that 1 implies $1^\p$.

	Now we assume $1^\p$ and proceed to prove 1. First we prove by induction
	that $H_1\ldots H_i$ is a subgroup of $G$ for each $i$. The assertion is
	true for $i=1$. For larger $i$, assume $H_1\ldots H_{i-1}$ is a subgroup
	and compute
		\begin{align*}
			(h_1\ldots h_i)(h_1^\p \ldots h_i^\p) &= (h_1\ldots h_{i-1})\Big(h_i(h_1^\p \ldots h_{i-1}^\p) h_i^{-1}\Big) (h_ih_i^\p)\\
			&=(h_1\ldots h_{i-1})~^{h_i}(h_1^\p \ldots h_{i-1}^\p)(h_i h_i^\p) \\
			&=(h_1\ldots h_{i-1})~^{h_i}h_1^\p~^{h_i}h_2^\p \ldots~^{h_i}h_{i-1}^\p(h_i h_i^\p) \\
			&\in (H_1\ldots H_{i-1})~(H_1)(H_1H_2)\ldots(H_1\ldots H_{i-1})~H_i H_i \tag{By $1^\p$}\\
			&\subseteq (H_1\ldots H_{i-1})(H_1\ldots H_{i-1})\ldots(H_1\ldots H_{i-1})~H_i H_i \\
			&= (H_1\ldots H_{i-1}) H_i \tag{By ind. hyp.}
		\end{align*}
	showing that $H_1\ldots H_i$ is closed under multiplication.  Next the computation
		\begin{align*}
			(h_1\ldots h_i)^{-1} &= h_i^{-1}(h_1\ldots h_{i-1})^{-1}\\
			& = \Big(h_i^{-1}(h_1\ldots h_{i-1})^{-1}h_i\Big) h_i^{-1}\\
			&= \Big(~^{h_i^{-1}}(h_1\ldots h_{i-1})\Big)^{-1} h_i^{-1} \\
			&\in (H_1\ldots H_{i-1})^{-1} H_i  \tag{By $1^\p$} \\
			&= (H_1\ldots H_{i-1}) H_i \tag{By ind. hyp.}
		\end{align*}
	shows that $H_1\ldots H_i$ is closed under taking inverses as well, so it is a subgroup.
	Finally, to see that it is normal, first take $h_j\in H_j$ with $j > i$ and compute
		\begin{align*}
			~^{h_j}(h_1\ldots h_i) &= ~^{h_j}h_1~^{h_j}h_2\ldots~^{h_j}h_i \\
			&\in (H_1)(H_1H_2)\ldots(H_1\ldots H_i) \tag{By $1^\p$}\\
			&\subseteq(H_1\ldots H_i)(H_1\ldots H_i)\ldots(H_1\ldots H_i) \\
			&= (H_1\ldots H_i) \tag{subgp.}
		\end{align*}
	showing
		\begin{equation*} 
			~^{H_j}(H_1\ldots H_i) \subseteq H_1\ldots H_i \tag{$*$}
		\end{equation*}
	which also holds for $j \le i$ since in that case $H_j \subseteq
	H_1\ldots H_i$ is an inclusion of subgroups. Then for arbitrary $g\in G$ the
	relation
		$$~^g(H_1\ldots H_i) \subseteq H_1\ldots H_i$$
	is seen to hold by writing $g$ as a product of factors of the form
	$h_j\in H_j$ for various $j$---such a factorization exists since the
	$H_i$ are assumed to generate $G$---and then applying $(*)$ repeatedly
	for each factor. 
}

\bigskip

We turn to discuss the relationship between the notion of internal
$r$-SDP and the usual notion of internal semidirect product (see also
\cite{CarrascoCegarra}). First, saying that $G$ is an internal 2-SDP
of two of its subgroups $H_1,H_2$ is equivalent to saying that it is an
internal semidirect product of $H_1$ and $H_2$ in the usual sense.
Correspondingly, the notation of Definition \ref{iSDP} reduces to the
usual symbols used to denote this situation, shown here.
	$$G = H_1\rtimes H_2$$

In the cases $r > 2$, it follows from the definition of $r$-SDP 
that $H_1\ldots H_i$ is a 2-SDP of $H_1\ldots H_{i-1}$ and
$H_i$ for $i > 1$. Using the above notation in an obvious manner, one
has for any $j$
	\begin{align*}
		H_1\ldots H_j &= \big(H_1\ldots H_{j-1}\big)\rtimes H_j \\
		&= \Big(\big(H_1\ldots H_{j-2}\big)\rtimes H_{j-1}\Big)\rtimes H_j \\
		&\hspace{6.5pt}\vdots \\
		&= \bigg(\ldots\Big(\big(H_1\rtimes H_2\big)\rtimes H_3 \Big) \rtimes
			\ldots \rtimes H_{j-1}\bigg)\rtimes H_j
	\end{align*}
and the case $j = r$ appears as follows.
	\begin{equation*}
		G = \bigg(\ldots\Big(\big(H_1\rtimes H_2\big)\rtimes H_3 \Big) \rtimes
			\ldots \rtimes H_{r-1}\bigg)\rtimes H_r \tag{$**$}
	\end{equation*}
A group $G$ of the form $(**)$ may be called an {\it iterated $2$-SDP}, 
and the above argument shows that every $r$-SDP is an
iterated 2-SDP. The converse statement is false, however, as the
subgroup $H_1\ldots H_i$ of an iterated 2-SDP is not required to be
normal in $G$ but only in $H_1\ldots H_{i+1}$. 

\bigskip

The preceding discussion leads to one last formulation, in terms of
normal towers, of the notion of internal $r$-SDP. It is helpful to do this first
for the more general notion of iterated 2-SDP. Say a normal tower
in $G$ 
	$$\xymatrix@C=0pt{\{1\} & =&  N_0 & \lhd & N_1 \ar@{->>}[d] & \lhd & N_2 \ar@{->>}[d] 
		& \lhd & \ldots & \lhd & N_r\ar@{->>}[d] & = & G \\
					&&&& N_1/N_0 && N_2/N_1 &&&& N_r/N_{r-1}}$$
is {\it completely split} if each of the short exact sequences
	\begin{equation*}\label{towersplitting}
		\xymatrix{N_{i-1} \ar@{^(->}[r] &  N_i \ar@{->>}[r] & N_i / N_{i-1}}
	\end{equation*}
is split.  Giving a decomposition $(**)$ of $G$ as an iterated 2-SDP 
is equivalent to giving a completely split normal tower in $G$ together
with a choice of splittings of each of the short exact sequences above. 
Indeed, such a tower is obtained from $(**)$ by setting
	$$N_i := H_1\ldots H_i$$
and taking for splittings the following canonical inclusions.
	$$\xymatrix{N_i = H_1\ldots H_i & ~H_i \cong N_i / N_{i-1} \ar@{_(->}[l]}$$
Conversely, given such a tower and a choice of splittings $\xymatrix{N_i
& ~N_i / N_{i-1}\ar@{_(->}[l]}$, the images of these splittings
constitute a family of subgroups $H_i$ of $G$ such that there is an
iterated 2-SDP decomposition $(**)$.

The notion of internal $r$-SDP is now recovered by imposing on a
completely split normal tower with chosen splittings the additional
requirement that each term of the tower is normal not only in the next
term, but indeed in $G$.

\bigskip

\rmk{
	This final remark discusses { internal} versus { external} $r$-SDPs. 
	First recall the following facts about the usual 2-fold notions.
	The notation $G = H_1\rtimes H_2$ for an {\it internal} semidirect
	product decomposition is only unambiguous in a context in which the
	group $G$ is given and the groups $H_1,H_2$ are specified as subgroups
	of $G$. Otherwise, the right-hand side generally does not uniquely
	determine the left-hand side, that is, there may be distinct
	nonisomorphic groups $G$ which are internally semidirect products of
	copies of $H_1$ and $H_2$. Indeed, it is well-known that $G$ is
	only determined by $H_1$ and $H_2$ together with a third piece of data,
	namely an action of $H_2$ on $H_1$ by automorphisms. Given such an
	action, the group $G$ may be reconstructed from $H_1$ and $H_2$ via the
	familiar {\it external semidirect product} construction.
\newline\indent 
	Quite analogous statements hold for $r$-SDPs (as well as iterated 2-SDPs). 
	The notation $G = H_1\rtimes\ldots\rtimes H_r$ indicates only
	that a given group $G$ stands in a certain relation to certain of its subgroups
	$H_i$---namely, the conditions of Definition \ref{iSDP} hold. In
	the subsequent sections of this article, an {\it external
	$r$-semidirect product} construction is presented yielding
	a well-defined group $G$ from abstractly given groups $H_i$ and certain
	extra data, such that $G$ is an internal SDP of the $H_i$ and such that
	(up to isomorphism) all such $G$ are obtained in this manner. 
}

	\section{Nonassociative External $r$-fold SDPs}				

	\label{ExternalSDPsHsPhis}

In \cite{CarrascoCegarra}, the authors analyze a general internal SDP
	$$G = H_1 \rtimes \ldots\rtimes H_r$$ 
and extract from it certain data describing the interaction of the
groups $H_i$ under the group law of $G$. They dub such a collection of
data a {\it total system}.

In order to find axioms characterizing such total systems, we take the
approach of assuming abstractly given groups $H_i$ and an abstractly
given total system, to be assembled into an {\it external $r$-SDP}
via a construction generalizing the usual external 2-SDP construction.
As also happens in the $r=2$ case, for general $r$ certain conditions
must be placed on a total system in order to guarantee associativity of
the resulting multiplication law. We give a method for generating these associativity conditions
algorithmically in section \ref{BracketAxiomsPhis}.

In preparation, we first work with a weaker nonassociative
construction, using the data in the total system to give to the set
	$$G := H_1 \times \ldots \times H_r$$
the structure of a {\it unital magma}, that is, a set with a
multiplication $\mu$ having a two-sided identity but otherwise not
required to satisfy any particular axioms. 


If the total system is such that the multiplication $\mu$ is fully
associative, it will then follow immediately from the definition of
$\mu$ and the assumption that the $H_i$ are groups that $G$ is a group
under $\mu$. However, the construction of $\mu$ is carried out without
making use of the existence of inverses in the $H_i$, and so the results
of this chapter also make sense for monoids $H_i$. The
construction can be further generalized to the case of magmas
$H_i$---even nonunital ones---with some straightforward adjustments, but
for expedience the $H_i$ are simply assumed to be groups throughout.

\bigskip

Let abstract groups $H_i$ be given for $1 \le i \le r$.
The following definition is due to Carrasco and Cegarra, introduced in
\cite{CarrascoCegarra}.

\defn{\label{totalsystem}
	An {\it $r$-total system}
		$$S^{(r)} = \Big\{ ~ \left\{ \phi_k^j \right\}_{r \ge k > j \ge 1} ~ ; ~
			\left\{~[\cdot,\cdot]_{kj}^i~\right\}_{r \ge k > j > i \ge 1} ~ \Big\}$$
	consists of a family of maps
		$$\phi_k^j : H_k \longrightarrow \mathrm{Maps}(H_j,H_j), ~~~ r\ge k > j \ge 1$$
		$$ h_k \mapsto \big(h_j \mapsto {}^{\phi_k^j(h_k)}h_j\big)$$
	where $\mathrm{Maps}(H_j,H_j)$ is the set of all set-theoretic self-maps of
	$H_j$, and a family of brackets
		$$[\cdot,\cdot]_{kj}^i : H_k \times H_j \longrightarrow H_i, ~~~ r \ge k > j > i \ge 1.$$
}

\rmk{
	The idea behind this definition is that the structure of an SDP of the $H_i$ should
	be determined via relations of the form
		\begin{equation*}
			h_k h_j = [h_k,h_j]_{kj}^1~[h_k,h_j]_{kj}^2~\ldots~[h_k,h_j]_{kj}^{j-1}
				\big({}^{\phi_k^j(h_k)}h_j\big)h_k
		\end{equation*}
	with $h_j\in H_j$ and $h_k\in H_k$ where $k > j$.  
	The reader is invited to compare condition $1^\p$ of Proposition \ref{checktermbyterm}.
\newline\indent
	Once the conditions on
	$S^{(r)}$ equivalent to associativity of $G$ have been determined, it is
	a simple matter to show that, under those conditions, $G$ is isomorphic
	to the group $U$ obtained as the quotient of the (group-valued) free
	product $H_1*\ldots *H_r$ by the above relations. 
\newline\indent
	It is tempting to try to start with $U$ defined in this way as a candidate
	for external SDP, but from
	there it is not easy to determine what algebraic conditions should
	be imposed on the total system $S^{(r)}$ so that the natural
	homomorphisms $H_i\hra U$ are injective. The difficulty seems to be that
	the right conditions are precisely the associativity conditions---but 
	in $U$ associativity is already artificially imposed (via the notion
	of free product of groups), leading to confusion.  Instead, we proceed
	by artificially imposing embeddedness of the $H_i$ and initially foregoing associativity.
}

The following normalization conditions will guarantee, among other
things, that the multiplication $\mu$ (still to be constructed)
restricts on $H_i$ to its originally given group law for each $i$.

\defn{
	Say the $r$-total system $S^{(r)}$ is {\it normalized} if
	for all $r \ge k > j > i \ge 1$ the following requirements are satisfied.
	\begin{itemize}
		\item[i)]{$\phi_j^i(1)$ is the identity map of $H_i$.}
		\item[ii)]{For all $h_j\in H_j$, $\phi_j^i(h_j)$ fixes the identity element $1\in H_i$.}
		\item[iii)]{For all $h_j\in H_j$, $[1,h_j]_{kj}^i = 1 \in H_i$.}
		\item[iv)]{For all $h_k\in H_k$, $[h_k,1]_{kj}^i = 1 \in H_i$.}
	\end{itemize}
}


\bigskip

Now assume given, in addition to the groups $H_i$, a normalized total
system $S^{(r)}$. In order to use this data to define the desired multiplication
$\mu$ on the product set
	$$G := H_1 \times \ldots \times H_r$$
the following preliminary definitions are needed.

\bigskip

\defn{ 
	The {\it rank} of an $r$-tuple $(h_1,\ldots,h_r)\in G$ is the largest
	index $j$ for which $h_j$ is nontrivial. Similarly, its {\it corank} is
	the smallest index $j$ for which $h_j$ is nontrivial.  Use the symbols $R_j$ and
	$R_j^\prime$ to denote the sets
			$$R_j := \Big\{~ (h_1,\ldots,h_j,1,\ldots,1) ~\Big\vert~ h_i \in H_i, ~1\le i \le j ~\Big\}$$
			$$R_j^\prime := \Big\{~ (1,\ldots,1,h_j,\ldots,h_r) ~\Big\vert~ h_i \in H_i, ~j \le i \le r ~\Big\}$$
	of all $r$-tuples of rank $\le j$, respectively of all $r$-tuples of corank $\ge j$. 
}


\defn{
	Let $A^{(2)}$ denote the subset of $G^{\times 2} = G\times G$ consisting of 
	pairs of $r$-tuples of the form
		$$(h_1, \ldots, h_i, 1, \ldots, 1), (1, \ldots, 1, h_i^\p, \ldots, h_r^\p)$$
	for any $i$ with $1\le i \le r$, which will be referred to as {\it noninterfering pairs}. That is,
		$$A^{(2)} := \bigcup_{i=1}^r R_i \times R_i^\prime.$$
	More generally, define $A^{(k)}$ to be the set of {\it noninterfering $k$-tuples}
	of $r$-tuples in $G$, that is, the subset of $G^{\times k}$ in which each component
	$r$-tuple has rank less than or equal to the corank of the next.
	Finally, let $A$ denote the set of all noninterfering families of $r$-tuples, 
	that is, $A$ is the disjoint union of the $A^{(k)}$ for all $k\ge 2$.
}


A magma multiplication $\mu$ will now be defined on the product set $G$
in several steps. In the initial step, a multiplication $\mu_A$ is
defined on $A^{(2)}$, extending immediately to $A$ in a fully
associative manner. Then the existence of $\mu_A$ enables us to
introduce some indispensible notational conventions. Finally,
multiplications $\mu_j$ on $R_j\times R_j$ are defined inductively, and
the final $\mu_r$ will be the desired multiplication $\mu$. It will not
be obvious initially how these various multiplications are related, but
as a consequence of the normalization conditions on $T^{(r)}$, they will
be shown to be pairwise coextensive. In particular, $\mu = \mu_r$ will
extend all of them to the domain $G\times G$.

\bigskip

\defn{
	Define the
	multiplication $\mu_A$ on $A^{(2)}$ via the following formula.
	\begin{equation*} 
		\mu_A\Big((h_1, \ldots, h_i, 1, \ldots, 1),(1, \ldots, 1, h_i^\p, \ldots, h_r^\p)\Big) := 
		(h_1, \ldots, h_ih_i^\p, \ldots, h_r^\p)
	\end{equation*}
}

\bigskip

\begin{prp} \label{propmuA}
	The following are immediate from the definition of $\mu_A$.
	\begin{enumerate}
		\item{The $r$-tuple $(1,\ldots,1)$ is a two-sided unit for $\mu_A$.}
		\item{$\mu_A$ restricts on each $H_i$ to the originally given multiplication law for $H_i$.}
		\item{$\mu_A$ is fully associative on all families of $r$-tuples in $A$.}
	\end{enumerate}
\end{prp}

\bigskip

\rmk{
	In the second assertion of this proposition, the $H_i$ are regarded as embedded in $G$
	via the usual canonical inclusions $h_i \hra (1, \ldots, 1, h_i, 1, \ldots, 1)$.
	Also note the third assertion requires associativity of the $H_i$.
}

\bigskip

It is convenient to take advantage of $\mu_A$ to introduce the following notation. 
We shall usually denote {\it elementary $r$-tuples}, that
is, $r$-tuples with a single nontrivial entry, by the abbreviation
	$$h_i := (1, \ldots, h_i, \ldots, 1)$$
for $h_i\in H_i$. In using this abbreviation, the subscript
indicates which entry in the elementary $r$-tuple is supposed to be
nontrivial. 
	
Next, the multiplication $\mu_A$ will be denoted with a dot, and so a dot is used 
when multiplying elements from different groups $H_i, H_j$ with $i <
j$, as in $h_i\cdot h_j$. When multiplying elements from the same group
$H_i$, the relation
	$$h_i \cdot h_i^\p = (1, \ldots, h_ih_i^\p, \ldots, 1) = h_ih_i^\p$$
holds, where the multiplication in $H_i$ is indicated by juxtaposition.

Using this alternate notation, 
every $r$-tuple can be written as a fully associative
$\mu_A$-product of elementary $r$-tuples with increasing indices
	$$(h_1,\ldots,h_r) = h_1 \cdot\ldots\cdot h_r$$
and an $r$-tuple in $R_j$ takes the convenient form
	$$(h_1,\ldots,h_j, 1, \ldots, 1) = h_1 \cdot\ldots\cdot h_j.$$

When it is desirable to be explicit concerning units, the
unit element of the group $H_i$ will be denoted by $1_i$ and the
trivial $r$-tuple $(1_1,\ldots,1_r)$ by $1^{(r)}$. The perhaps
strange-looking equations
	$$1_i = 1^{(r)} = 1_j$$
holding for all $i,j$ are a result of these notational conventions.  When it
is desirable to be less explicit, the symbol $1$ will be used to denote any of
these.

\bigskip

In order to proceed with the construction of $\mu$, the data
in the total system $S^{(r)}$ will be encoded in the following operators.

\defn{\label{conjcommop}
	For $r \ge k > j \ge 1$, define
	the {\it conjugation} operator
		$$\phi_{(\cdot)}(\cdot): H_k\times H_j \longrightarrow R_j$$
	\begin{equation*} \label{foperator}
		\phi_{h_k}(h_j) := [h_k,h_j]_{kj}^1\cdot[h_k,h_j]_{kj}^2\cdot\ldots
			\cdot[h_k,h_j]_{kj}^{j-1}\cdot {}^{\phi_k^j(h_k)}h_j
	\end{equation*}
	and for later use also define the following {\it commutator bracket}.
		$$[h_k,h_j] := \phi_{h_k}(h_j)\cdot h_j^{-1}$$
		$$ = [h_k,h_j]_{kj}^1\cdot\ldots
			\cdot[h_k,h_j]_{kj}^{j-1}\cdot [h_k,h_j]_{kj}^j$$
	where the following convenient notation has been introduced.
		$$[h_k,h_j]_{kj}^j := {}^{\phi_k^j(h_k)}h_jh_j^{-1}$$
}

\rmk{
	Note the $r$-tuples $\phi_{h_k}(h_j)$ and $[h_k,h_j]$ belong to $R_j$. Of course
	$\phi$ depends on $k$ and $j$, but the particular $\phi$ being used will be inferred
	from the indices of the arguments $h_k, h_j$. The normalization
	conditions for $S^{(r)}$ translate into the conditions $$\phi_{1_k}(h_j) =
	h_j$$
		$$\phi_{h_k}(1_j) = 1^{(r)}$$
	for all $h_j\in H_j$, $h_k\in H_k$.
}

\defn{ \label{BigMuDef}
	Now inductively define multiplications $\mu_k$ on $R_k$, starting 
	with the observation that there is already a well-defined
	multiplication 
		$$\mu_1: R_1\times R_1\rightarrow R_1$$ 
	given by the group
	law for $H_1$.  Assuming there is a well-defined multiplication
	$\mu_k: R_k\times R_k\rightarrow R_k$ for some $k$ with $1 \le k \le r-1$, define a
	multiplication 
		$$\mu_{k+1}: R_{k+1}\times R_{k+1}\rightarrow R_{k+1}$$ 
	by the formula
	\begin{eqnarray*} 
		\mu_{k+1}\Big((a_1\cdot\ldots\cdot a_{k+1}), (b_1\cdot\ldots\cdot b_{k+1})\Big) := ~~~~~~~~~~~~~~ 
			\nonumber\\
		\mu_k\Bigg((a_1\cdot\ldots\cdot a_k), ~
			\mu_k\Big( \phi_{a_{k+1}}(b_1),  \ldots, \phi_{a_{k+1}}(b_k) \Big)\Bigg) 
				\cdot(a_{k+1}b_{k+1})
	\end{eqnarray*}
	where the following recursive definition
	\begin{eqnarray*} 
		\mu_i(w_1, \ldots, w_j) := \mu_i\Big(\mu_i(w_1, \ldots, w_{j-1}), w_j\Big)\nonumber\\
		 = \mu_i\Big(\ldots\mu_i\big(\mu_i(w_1,w_2), w_3\big) \ldots, w_j\Big)
	\end{eqnarray*}
	is used for $w_1,\ldots,w_j\in R_i$.
	When $j$ equals 1, $\mu_i(w_1)$ is understood to be simply $w_1$. 
	The full multiplication law $\mu_r$ is obtained by taking $k+1=r$ in the
	above formula. 
}

\begin{lem} \label{muA}
	For each $k$ with $1\le k \le r$, 
	the multiplications $\mu_k$ and $\mu_A$ agree on their common domain
	$(R_k\times R_k)\cap A^{(2)}$.  In particular, $\mu_r$ is an extension of $\mu_A$.
\end{lem}
\prf{Proof.}{
	Proceed by induction.  To begin, note that $\mu_1$ is just the group law for the
	embedded copy of $H_1$, and so it agrees with $\mu_A$ on $R_1 = H_1$.

	Now assume for induction that $\mu_k$ agrees with $\mu_A$ on 
	$(R_k\times R_k)\cap A^{(2)}$.  One proceeds as follows first in the case $1\le i < k+1$.
		\begin{align*}
			&\mu_{k+1}\Big((a_1\cdot\ldots\cdot a_i), (b_i\cdot\ldots\cdot b_{k+1})\Big) = \\
			&= \mu_{k+1}\Big((a_1\cdot\ldots\cdot a_i\cdot 1_{k+1}), (b_i\cdot\ldots\cdot b_{k+1})\Big)\\
			&= \mu_k\Bigg((a_1\cdot\ldots\cdot a_i),
				\mu_k\Big(\phi_1(1_1),\ldots,\phi_1(1_{i-1}), \phi_1(b_i), \ldots, \phi_1(b_k) \Big)\Bigg)\cdot(1_{k+1}b_{k+1}) \\
			&= \mu_k\Big((a_1\cdot\ldots\cdot a_i), 
				\mu_k\left(1^{(r)},\ldots,1^{(r)}, b_i, \ldots, b_k\right)\Big)\cdot b_{k+1}
		\end{align*}
	Now in the term $\mu_k\left(1^{(r)},\ldots,1^{(r)}, b_i, \ldots,
	b_k\right)$, the arguments constitute a noninterfering
	family belonging to $A$.
	Hence by the inductive hypothesis, the result of evaluating $\mu_k$ 
	on the corresponding parenthesized word
	coincides with the result of evaluating $\mu_A$ on the family.
	This result is simply $b_i\cdot\ldots\cdot b_k$,
	so the above chain of equalities can be continued with
		\begin{align*}
			\mu_k\Big((a_1\cdot\ldots\cdot a_i), (b_i\cdot\ldots\cdot b_k)\Big)\cdot b_{k+1} &= 
				\mu_A\Big((a_1\cdot\ldots\cdot a_i),	(b_i\cdot\ldots\cdot b_k)\Big)\cdot b_{k+1} \\
			&= a_1\cdot\ldots\cdot a_{i-1} \cdot a_ib_i\cdot b_{i+1}\cdot\ldots\cdot b_{k+1} \\
			&= \mu_A\Big((a_1\cdot\ldots\cdot a_i), (b_i\cdot\ldots\cdot b_{k+1})\Big)
		\end{align*}
	as desired.  

	For the case $i = k+1$ one also obtains
		\begin{align*}
			&\mu_{k+1}\Big((a_1\cdot\ldots\cdot a_{k+1}), b_{k+1}\Big) =  \\
		 	&= \mu_k\Bigg((a_1\cdot\ldots\cdot a_k), \mu_k\Big(\phi_{a_{k+1}}(1_1),\ldots,\phi_{a_{k+1}}(1_k)\Big)\Bigg)\cdot(a_{k+1}b_{k+1}) \\
			&= \mu_k\Big((a_1\cdot\ldots\cdot a_k), 
				\mu_k\left(1^{(r)},\ldots,1^{(r)}\right)\Big)\cdot(a_{k+1}b_{k+1}) \\
			&= \mu_k\Big((a_1\cdot\ldots\cdot a_k), 
				\mu_A\left(1^{(r)},\ldots,1^{(r)}\right)\Big)\cdot(a_{k+1}b_{k+1}) \\
			&= \mu_k\Big((a_1\cdot\ldots\cdot a_k), 
				1^{(r)}\Big)\cdot(a_{k+1}b_{k+1}) \\
			&= \mu_A\Big((a_1\cdot\ldots\cdot a_k), 
				1^{(r)}\Big)\cdot(a_{k+1}b_{k+1}) \\
			&= a_1\cdot\ldots\cdot a_k\cdot a_{k+1}b_{k+1} \\
			&= \mu_A\Big((a_1\cdot\ldots\cdot a_{k+1}), b_{k+1}\Big)
		\end{align*}
	as desired.  
}

\begin{cor}
	The multiplications $\mu_k$ for $1 \le k \le r$ form an increasing sequence of extensions.
\end{cor}
\prf{Proof.}{
	Fixing $k\ge 1$, proceed as follows to show that $\mu_{k+1}$ extends $\mu_k$.
		\begin{align*}
			&\mu_{k+1}\Big((a_1\cdot\ldots\cdot a_{k}), (b_1\cdot\ldots\cdot b_{k})\Big) = \\
			&= \mu_{k+1}\Big((a_1\cdot\ldots\cdot a_{k}\cdot 1_{k+1}), (b_1\cdot\ldots\cdot b_{k}\cdot 1_{k+1})\Big) \\
			&= \mu_{k}\Bigg((a_1\cdot\ldots\cdot a_{k}), 
				\mu_{k}\Big( \phi_1(b_1), \phi_1(b_2), \ldots, \phi_1(b_{k}) \Big)\Bigg)\cdot(1_{k+1}1_{k+1}) \\
			&= \mu_{k}\Bigg((a_1\cdot\ldots\cdot a_{k}), 
				\mu_{k}\Big( b_1, \ldots, b_{k} \Big)\Bigg) 
		\end{align*}
	Now observe that the arguments of the term $\mu_{k}(b_1, \ldots, b_{k})$ constitute
	a noninterfering family	belonging to $A$. Hence
	by the previous lemma (Lemma \ref{muA}), this term coincides with 
		$$\mu_A\big( b_1, \ldots, b_{k} \big) = b_1\cdot\ldots\cdot b_{k}$$
	and consequently $\mu_{k+1}$ restricts to $\mu_{k}$ on $R_{k}\times R_{k}$, as desired.
}

\begin{cor} \label{unit}
	The following hold for the multiplication $\mu =\mu_r :G\times G\ra G$.
	\begin{enumerate}
		\item{The $r$-tuple $1^{(r)} = (1,\ldots,1)$ is a two-sided unit for $\mu$.}
		\item{The multiplication $\mu$ restricts on the embedded copy of $H_i$
			to the originally given group law for $H_i$.}
	\end{enumerate}
\end{cor}
\prf{Proof.}{
	Since $\mu$ extends $\mu_A$ by Lemma \ref{muA}, the assertions follow from
	the corresponding assertions for $\mu_A$ (see Proposition \ref{propmuA}).
%
}

\bigskip

	\section{Conjugation Operators and Associativity Relations}		
	\label{SomeAssocsPhis}

Let groups $H_1, \ldots, H_r$ be given along with a normalized total
system $S^{(r)}$, and let $G = H_1\times\ldots\times H_r$ have the
multiplication $\mu$ constructed from these data as in the previous
section. Also recall the partially defined multiplications $\mu_A$ and
$\mu_k$, which are coextensive with $\mu$ as proved in the previous
section.

In the present section, we first extend the conjugation operators
$\phi_{a_k}$ in a natural manner suggested by the recursive definition
of $\mu$ (Definition \ref{BigMuDef}). This leads to additional
simplifications of the notation. Then we show how associativity
relations for $\mu$ can be translated into algebraic conditions on the
conjugation operators $\phi_{a_k}$. Finally, two simple associativity
relations are recorded that hold for general $\mu$, that is, without any
special conditions on the total system $S^{(r)}$. These are used later
in the determination of necessary and sufficient conditions for full
associativity of $\mu$. In fact, there are many more such generally
valid associativities, and we mention some of these but do not pursue
the matter further.

\bigskip

Since the multiplications $\mu_A$ and $\mu_k$ are compatible, the dot
notation, previously used only for $\mu_A$, is hereafter extended to all
$\mu_k$ without danger of ambiguity. Of course the formula for $\mu_A$
still only applies on the domain $A$. Moreover, since the $\mu_k$ are
generally not associative, we will have to keep track of parentheses. 

\bigskip

The operators $\phi_{a_k}$ may now also be extended as follows.

\bigskip

\defn{ \label{extendedf}
	For $a_k\in H_k$, define an operation 
		$$\phi_{a_k}:R_{k-1} \longrightarrow R_{k-1}$$
	simultaneously extending all of the operations
		$$\phi_{a_k}:H_i \longrightarrow R_i ~\mfor~ 1 \le i \le k-1$$
	via the following recursive formula.
		\begin{equation*} 
			\phi_{a_k}(b_1\cdot\ldots\cdot b_{k-1}) := 
				\phi_{a_k}(b_1\cdot\ldots\cdot b_{k-2})\cdot \phi_{a_k}(b_{k-1})
		\end{equation*}
		\begin{equation*}
			= \Bigg(\ldots\bigg(\Big(\phi_{a_k}(b_1)\cdot \phi_{a_k}(b_2)\Big)\cdot 
			\phi_{a_k}(b_3)\bigg)\cdot \ldots\cdot \phi_{a_k}(b_{k-2})\Bigg)\cdot \phi_{a_k}(b_{k-1})
		\end{equation*}
	This is indeed an $r$-tuple in $R_{k-1}$ by the previous definition of $\phi_{a_k}$.
}

\rmk{ \label{premutiny}
	Note that $\phi_{a_k}(b_1\cdot\ldots\cdot b_{k-1})$ as defined above is precisely the term
			$$\mu_{k-1}\Big(\phi_{a_k}(b_1), \ldots, \phi_{a_k}(b_{k-1})\Big)$$
	appearing in Definition \ref{BigMuDef}.
	Indeed, with this notational convention, the inductive definition of
	$\mu$ takes the simpler form
		\begin{equation*} 
			(a_1\cdot\ldots\cdot a_k) \cdot (b_1\cdot\ldots\cdot b_k) :=
				\Big((a_1\cdot\ldots\cdot a_{k-1}) \cdot
				\phi_{a_k}(b_1\cdot\ldots\cdot b_{k-1}) \Big) \cdot(a_kb_k)
		\end{equation*}
	or equivalently 
		\begin{equation} \label{premu2}
			(u\cdot a_k) \cdot (v\cdot b_k) := \Big(u \cdot \phi_{a_k}(v)\Big) \cdot(a_kb_k)
		\end{equation}	
	where 
		$$u := a_1\cdot\ldots\cdot a_{k-1}$$
		$$v := b_1\cdot\ldots\cdot b_{k-1}.$$
\indent	
	Finally, the following 
	useful identity is obtained from \eqref{premu2} 
	by taking $u$ and $b_k$ to be trivial.  It holds for all $v$ of rank $\le k-1$.
	\begin{equation} \label{tiny}
		a_k\cdot v = \phi_{a_k}(v) \cdot a_k
	\end{equation}
}

\defn{
	Let $U, V, W$ be subsets of $G$.  Let the symbol
		$${\bf A}[U,V,W]$$
	stand for the family of all associativity relations of the form
		$$(u\cdot v)\cdot w = u\cdot (v\cdot w)$$
	for $u\in U, v\in V, w\in W$.  
\newline\indent
	In the special case in which one of $U,V,W$ is $R_k$ or $H_k$, 
	the symbol $(k)$ or $k$ will be substituted for it respectively in the above notation.
	Here is	an example.
		$${\bf A}[k,(j),(i)] := {\bf A}[H_k, R_j, R_i]$$
	The conditions of the special form ${\bf A}[k,j,i]$ as will be referred to as
	\emph{elementary associativity relations}. 
	The following abbreviation may also be used for emphasis.
		$${\bf A}[[k]] := {\bf A}[R_k,R_k,R_k]$$
}

\rmk{
	Note one trivially has all associativity relations of the form ${\bf
	A}[i,j,k]$ for $i\le j \le k$ by the associativity of $\mu$ on the
	domain $A$ (see the previous section).  
\newline\indent
	More generally, it can be shown that ${\bf A}[R_i, R_j^\p, G]$ holds for
	any $i < j$. Other general associativities hold as well, such as ${\bf
	A}[R_k, R_k, R_k^\p]$ and ${\bf	A}[R_k, R_k^\p, R_k^\p]$ for any $k$.
	Since we will only need the two simple results
	proved in Proposition \ref{posspit} below, the ones mentioned in this remark are left
	as exercises.
}

\begin{lem} \label{translate} Associativity conditions for $\mu$ can
	be translated into conditions on the operators $\phi_{a_k}$ in the
	following ways.
	\begin{enumerate}
		\item{For $V,W\subseteq R_{k-1}$, the condition ${\bf A}[H_k,V,W]$
			is equivalent to the assertion that for all $a_k \in H_k$, $v\in V$, $w\in W$ one has
				$$\phi_{a_k}(v\cdot w) = \phi_{a_k}(v) \cdot \phi_{a_k}(w).$$}
		\item{For $W\subseteq R_{k-1}$, the condition ${\bf A}[H_k,H_k,W]$
			is equivalent to the assertion that for all
			$a_k, b_k\in H_k$ and $w\in W$ one has
				$$\phi_{a_kb_k}(w) = \big(\phi_{a_k}\circ \phi_{b_k}\big)(w).$$}
	\end{enumerate}
\end{lem}
\prf{Proof.}{
	To prove the first assertion, assume the associativity condition
		$$(a_k \cdot v)\cdot w = a_k \cdot (v\cdot w)$$
	and argue as follows.
		\begin{align*}
			\phi_{a_k}(v\cdot w)\cdot a_k &= a_k\cdot(v\cdot w) \phantom{\Big(}\tag{By \eqref{tiny}}\\
			&= (a_k \cdot v)\cdot w \phantom{\Big(} \\
			&= \Big(\phi_{a_k}(v)\cdot a_k\Big) \cdot w \tag{By \eqref{tiny}} \\
			&= \Big(\phi_{a_k}(v)\cdot \phi_{a_k}(w)\Big) \cdot a_k \tag{By \eqref{premu2}}
		\end{align*}
	Now the first $k-1$ entries of the far left-hand side are the nontrivial
	entries of \mbox{$\phi_{a_k}(v\cdot w)$}, and the first $k-1$ entries of the far
	right-hand side are the nontrivial entries of $\phi_{a_k}(v)\cdot \phi_{a_k}(w)$.
	So one concludes
		$$\phi_{a_k}(v\cdot w) = \phi_{a_k}(v) \cdot \phi_{a_k}(w)$$
	as claimed.  The proof of the converse is similar but easier.
	
	The proofs of the two directions of the second assertion are quite similar to those of the first,
	so we just give the relevant calculation and leave the rest to the reader.
		\begin{align*}
			\phi_{a_kb_k}(w)\cdot (a_kb_k) &= (a_kb_k)\cdot w \phantom{\Big(}\tag{By \eqref{tiny}} \\
			&= a_k\cdot (b_k \cdot w) \tag{By hypothesis} \phantom{\Big(}\\
			&= a_k \cdot \Big(\phi_{b_k}(w)\cdot b_k\Big) \tag{By \eqref{tiny}} \\
			&= \big(\phi_{a_k}\circ \phi_{b_k}\big)(w) \cdot (a_kb_k)\phantom{\Big(}\tag{By \eqref{premu2}}
		\end{align*}
}


The following simple proposition, for use in the next section, gives
just two of the many associativity relations holding for $\mu$ in
general, that is, without any special conditions on the total system
$S^{(r)}$.

\begin{prp} \label{posspit}
	The following associativity conditions hold in general.
	\begin{enumerate}
		\item{${\bf A}[(k-1),(k-1),k \hspace{1.2pt}]$ holds for $2 \le k \le r$.}
		\item{${\bf A}[k,(i),j\hspace{1pt}]$ holds for $1 \le i < j < k \le r$.}
	\end{enumerate}
\end{prp}
\prf{Proof.}{
	For the first assertion, compute as follows for $u,v\in R_{k-1}$ and $c_k\in H_k$.
		\begin{align*}
			u\cdot(v\cdot c_k) &= (u\cdot 1_k) \cdot (v\cdot c_k) \\
			&= \Big(u\cdot \phi_{1_k}(v)\Big) \cdot(1_kc_k) \tag{By \eqref{premu2}}\\
			&= (u\cdot v)\cdot c_k
		\end{align*}

	For the second assertion, it suffices by Lemma \ref{translate} to prove the identity
	\begin{equation*} 
		\phi_{a_k}(v\cdot c_j) = \phi_{a_k}(v)\cdot \phi_{a_k}(c_j)
	\end{equation*}
	for any $a_k\in H_k$, $v\in R_i$ and $c_j\in H_j$.
	This latter identity is an immediate consequence of the formula defining
	the extended conjugation operator $\phi_{a_k}$ in Definition \ref{extendedf}.
}


	\section{A Generating Family of Associativity Relations}		

	\label{GenAssocRelsNewPhis}

Once again let groups $H_1,\ldots, H_r$ be given along with a normalized
$r$-total system $S^{(r)}$, and let $G$ with multplication $\mu$ be
the corresponding external SDP as constructed in section \ref{ExternalSDPsHsPhis}. The
goal of this section is to prove the following two results.

\begin{thm} \label{mainlemma}
	Let $k$ with $r \ge k \ge 1$ be fixed.  In order for $\mu$ to be associative on $R_k$, 
	it suffices to require the following associativity conditions.
	\begin{enumerate}
		\item{${\bf A}[[k-1]]$, i.e., associativity of $\mu$ on $R_{k-1}$.}
		\item{${\bf A}[k,j,i]$ for all $i,j$ with $k \ge j \ge i \ge 1$.}
	\end{enumerate}
\end{thm}

From here a simple induction argument yields the following as a corollary.

\begin{thm}\label{mainassocthm}
	In order for $\mu$ to be associative on all of $G$, it suffices to assume only that the 
	elementary associativity conditions ${\bf A}[k,j,i]$ hold for all $i,j,k$ with
	$r \ge k \ge j \ge i \ge 1$.
\end{thm}

The proof of Theorem \ref{mainlemma} is carried out below after some preliminary work.

\bigskip

\rmk{	\label{synopsis}
	For the reader's convenience, the following
	synopsis gathers together the three basic algebraic identities used in the proofs below. 
	All hold without any additional conditions on the total system $S^{(r)}$.
	\begin{align} 
		\tag{\ref{tiny}}
		a_k\cdot v &= \phi_{a_k}(v) \cdot a_k & 
			\forall& a_k \in H_k, \forall v\in R_{k-1} \nonumber\\
		\label{spit} 
		\phi_{a_k}(v\cdot c_j) &= \phi_{a_k}(v)\cdot \phi_{a_k}(c_j) & 
			\forall& a_k \in H_k, \forall v\in R_i, \forall c_j\in H_j \mfor i<j<k \\
		\label{newmu} 
		(u \cdot a_k) \cdot (v \cdot b_k) &= u\cdot \phi_{a_k}(v) \cdot a_kb_k & 
			\forall& a_k,b_k \in H_k, \forall u,v\in R_{k-1}
	\end{align}
	Identity \eqref{tiny} appeared in the previous section under Remark \ref{premutiny}
	and has kept its line number. Identity
	\eqref{spit} is equivalent to the second assertion of Proposition \ref{posspit},
	as mentioned in the proof.  Finally, \eqref{newmu} is obtained from \eqref{premu2} of the 
	previous section by invoking the first assertion of Proposition \ref{posspit}, 
	allowing us to omit the larger pair of parentheses without ambiguity.
}

\bigskip





The next proposition gives general conditions under which the operation 
	$$\phi_{a_k}: R_{k-1}\longrightarrow R_{k-1}$$
is a homomorphism with respect to $\mu$.

\begin{prp} \label{hyp1}
	Let $k$ with $r \ge k > 1$ be fixed.  Assume ${\bf A}[[k-1]]$, and also
	assume ${\bf A}[k,j,i]$ for all $i,j$ such that $k > j \ge i \ge 1$.
	Then the condition 
		$${\bf A}[k,(k-1),(k-1)]$$
	holds, or equivalently the formula
	\begin{equation*} 
		\phi_{a_k}(u\cdot v) = \phi_{a_k}(u)\cdot \phi_{a_k}(v)
	\end{equation*}
	is valid for any $a_k\in H_k$ and $u,v\in R_{k-1}$.
\end{prp}

The heart of the proof in contained in the following lemma.  

\begin{lem} \label{scratch}
	Fix $i,j,k$ with $r \ge k > j > i \ge 1$.  Under the assumptions 
	$${\bf A}[[k-1]] \textrm{~~and~~} {\bf A}[k,(j-1),(j-1)]$$
	the following implications hold.
	\begin{enumerate}
  		\item{$\left.\begin{array}{l}
			{\bf A}[k,j,(i-1)]\\[2 pt]
			{\bf A}[k,j,i] \\[2 pt] 
			\end{array}\right\}  \implies {\bf A}[k,j,(i)]$}
		\item{$\left.\begin{array}{l} 
			{\bf A}[k,j,(j-1)] \\[2 pt]
			{\bf A}[k,j,j] \\[2 pt] 
		\end{array}\right\}  \implies {\bf A}[k,(j),(j)]$}
	\end{enumerate}
\end{lem}
\prf{Proof.}{
	In the following calculations, 
	 implicit use of the associativity condition ${\bf A}[[k-1]]$ is made
	through appropriate omission of parentheses.

	To prove the first assertion, take $a_k \in H_k$, $b_j \in H_j$, $w\in R_i$ and write
		$$w = w^\prime \cdot c_i$$
	where $w^\prime \in R_{i-1}$ and $c_i \in H_i$ are determined by $w$.  Then compute as follows.
	(The reader will notice that 
	\mbox{${\bf A}[k,(i-1),(i)]$} is invoked.  It follows {\it a fortiori} from the
	assumption \mbox{${\bf A}[k,(j-1),(j-1)]$} since $i$ is less than $j$.)
	\begin{align*}
		\phi_{a_k}(b_j \cdot w) &= \phi_{a_k}\Big(\phi_{b_j}(w)\cdot b_j\Big) \tag{By \eqref{tiny}} \\
		&= \phi_{a_k}\Big(\phi_{b_j}(w)\Big)\cdot \phi_{a_k}(b_j) \tag{By \eqref{spit}} \\
		&= \phi_{a_k}\Big(\phi_{b_j}(w^\prime)\cdot \phi_{b_j}(c_i)\Big)\cdot \phi_{a_k}(b_j) \tag{By \eqref{spit}}\\
		&= \phi_{a_k}\Big(\phi_{b_j}(w^\prime)\Big)\cdot \phi_{a_k}\Big(\phi_{b_j}(c_i)\Big)\cdot 
			\phi_{a_k}(b_j) \tag{By ${\bf A}[k,(i-1),(i)]$} \\
		&= \phi_{a_k}\Big(\phi_{b_j}(w^\prime)\Big)\cdot \phi_{a_k}\Big(\phi_{b_j}(c_i)\cdot b_j\Big) \tag{By \eqref{spit}} \\
		&= \phi_{a_k}\Big(\phi_{b_j}(w^\prime)\Big)\cdot \phi_{a_k}(b_j\cdot c_i) \tag{By \eqref{tiny}} \\
		&= \phi_{a_k}\Big(\phi_{b_j}(w^\prime)\Big)\cdot \phi_{a_k}(b_j)\cdot \phi_{a_k}(c_i) \tag{By ${\bf A}[k,j,i]$} \\
		&= \phi_{a_k}\Big(\phi_{b_j}(w^\prime)\cdot b_j\Big)\cdot \phi_{a_k}(c_i) \phantom{\Big(}\tag{By \eqref{spit}} \\
		&= \phi_{a_k}(b_j\cdot w^\prime)\cdot \phi_{a_k}(c_i) \phantom{\Big(}\tag{By \eqref{tiny}} \\
		&= \phi_{a_k}(b_j)\cdot \phi_{a_k}(w^\prime)\cdot \phi_{a_k}(c_i) \phantom{\Big(}\tag{By ${\bf A}[k,j,(i-1)]$} \\
		&= \phi_{a_k}(b_j)\cdot \phi_{a_k}(w^\prime \cdot c_i)\phantom{\Big(}\tag{By \eqref{spit}} \\
		&= \phi_{a_k}(b_j)\cdot \phi_{a_k}(w) \phantom{\Big(}
	\end{align*}
	By Lemma \ref{translate}, this identity is equivalent to the
	associativity condition ${\bf A}[k,j,(i)]$, which is thus proved.

	To prove the second implication,  take $a_k \in H_k$, $v, w\in R_j$ and write
		$$v = v^\prime \cdot b_j ~~~~~~~~w = w^\prime \cdot c_j$$
	where $v^\prime, w^\prime \in R_{j-1}$ and $b_j,c_j \in H_j$ are determined by $v,w$.
	Then compute as follows.
	\begin{align*}
		\phi_{a_k}(v \cdot w) &= 
			\phi_{a_k}\Big(v^\prime \cdot \phi_{b_j}(w^\prime)\cdot b_jc_j\Big) \tag{By \eqref{newmu}} \\
		&= \phi_{a_k}\Big(v^\prime \cdot \phi_{b_j}(w^\prime)\Big)\cdot \phi_{a_k}(b_jc_j) \tag{By \eqref{spit}}\\
		&= \phi_{a_k}(v^\prime) \cdot \phi_{a_k}\Big(\phi_{b_j}(w^\prime)\Big)\cdot \phi_{a_k}(b_j)\cdot \phi_{a_k}(c_j) 
					\tag*{
						$\left(\begin{array}{c} \text{By {\bf A}}[k,(j-1),(j-1)] \\[3pt]
						\text{and {\bf A}}[k,j,j]\end{array}\right)$
					}\\
		&= \phi_{a_k}(v^\prime) \cdot \phi_{a_k}\Big(\phi_{b_j}(w^\prime)\cdot b_j\Big)
			\cdot \phi_{a_k}(c_j) \tag{By \eqref{spit}}\\
	 	&= \phi_{a_k}(v^\prime) \cdot \phi_{a_k}\Big(b_j\cdot w^\prime\Big) \cdot \phi_{a_k}(c_j) \tag{By \eqref{tiny}}\\
		&= \phi_{a_k}(v^\prime) \cdot \phi_{a_k}(b_j)\cdot \phi_{a_k}(w^\prime)
			\cdot \phi_{a_k}(c_j) \tag{By ${\bf A}[k,j,(j-1)]$} \phantom{\Big(}\\
		&= \phi_{a_k}(v^\prime\cdot b_j)\cdot \phi_{a_k}(w^\prime \cdot c_j) \tag{By \eqref{spit}}\phantom{\Big(}\\
		&= \phi_{a_k}(v)\cdot \phi_{a_k}(w)  \phantom{\Big(}
	\end{align*}
	By Lemma \ref{translate}, this identity is equivalent to the
	associativity condition ${\bf A}[k,(j),(j)]$, as desired.
}

\prf{Proof of Proposition \ref{hyp1}.}{
	It suffices to prove
	by induction on $j$ that ${\bf A}[k,(j),(j)]$ holds for
	all $j$ with $1 \le j < k$.

	To start the induction, note that the condition
		$${\bf A}[k,(1),(1)] = {\bf A}[k,1,1]$$
	is among the hypotheses of the proposition.

	Now assume for induction that ${\bf A}[k,(j-1),(j-1)]$ holds
	for some $j$ with $k > j \ge 2$.  By the first part of 
	Lemma \ref{scratch} and an easy induction on $i$, one may conclude that 
	${\bf A}[k,j,(j-1)]$ holds.  By the second part of Lemma
	\ref{scratch}, one may then conclude that ${\bf A}[k,(j),(j)]$ holds.
}

\begin{prp} \label{hyp3}
	Let $k$ with $r \ge k > 1$ be fixed. Assume the hypotheses of
	Proposition \ref{hyp1} and additionally the conditions ${\bf A}[k,k,i]$
	for $k > i \ge 1$. Then the condition 
		$${\bf A}[k,k,(k-1)]$$
	holds, or equivalently the formula
	\begin{equation*} 
		\phi_{a_kb_k} = \phi_{a_k} \circ \phi_{b_k}
	\end{equation*}
	is valid on $R_{k-1}$ for any $a_k,b_k\in H_k$.
\end{prp}
\prf{Proof.}{
	It suffices to prove by induction on $i$ that, under the given hypotheses, 
	the assertion ${\bf A}[k,k,(i)]$
	holds for $i$ with $k > i \ge 1$.  

	To start the induction, note that the condition
		$${\bf A}[k,k,(1)] = {\bf A}[k,k,1]$$
	is among the given hypotheses.
 
	Now assume for induction that the condition ${\bf A}[k,k,(i-1)]$
	holds for some $i$ with $k > i \ge 2$.  
	Take $a_k,b_k \in H_k$ and $w\in R_i$ and write
		$$w = w^\prime \cdot c_i$$
	where $w^\prime \in R_{i-1}$ and $c_i \in H_i$ are determined by $w$.  Then compute
	as follows.
	\begin{align*}
		\phi_{a_kb_k}(w) &= \phi_{a_kb_k}(w^\prime) \cdot \phi_{a_kb_k}(c_i) \tag{By \eqref{spit}}\\
		&= \big(\phi_{a_k}\circ \phi_{b_k}\big)(w^\prime) \cdot \big(\phi_{a_k}\circ \phi_{b_k}\big)(c_i) 		
			\tag*{$\left(\begin{array}{c} 
				\text{By {\bf A}}[k,k,(i-1)] \\[3pt] 
				\text{and {\bf A}}[k,k,i]
			\end{array}\right)$}\\
		&= \phi_{a_k}\Big(\phi_{b_k}(w^\prime) \cdot \phi_{b_k}(c_i)\Big)\tag{By Proposition \ref{hyp1}}\\
		&= \phi_{a_k}\Big(\phi_{b_k}\big(w^\prime \cdot c_i\big)\Big)\tag{By \eqref{spit}}\\
		&= \big(\phi_{a_k} \circ \phi_{b_k}\big)(w) \phantom{\Big(}
	\end{align*}
	By Lemma \ref{translate}, this identity is equivalent to the 
	associativity condition ${\bf A}[k,k,(i)]$ as desired.
}

\prf{Proof of Theorem \ref{mainlemma}.}{
	First observe that the hypotheses of the theorem include those of Propositions
	\ref{hyp1} and \ref{hyp3} and so they hold. Now take $u,v,w \in R_k$
	and write
	\begin{align*}
		u &= u^\prime\cdot a_k \\
		v &= v^\prime\cdot b_k \\
		w &= w^\prime\cdot c_k
	\end{align*}
	where $u^\prime, v^\prime, w^\prime \in R_{k-1}$ and $a_k, b_k, c_k \in H_k$ are
	determined by $u,v,w$ respectively.  Then compute as follows.  (The absence of
	parentheses on the second line below is justified by Proposition \ref{posspit}
	together with the given hypothesis ${\bf A}[[k-1]]$.)
	\begin{align*}
		(u\cdot v)\cdot w &= 
			\Big(u^\prime \cdot \phi_{a_k}(v^\prime) \cdot a_kb_k\Big) \cdot w \tag{By \eqref{newmu}}\\
		&= u^\prime \cdot \phi_{a_k}(v^\prime) \cdot \phi_{a_kb_k}(w^\prime) \cdot a_kb_kc_k \phantom{\Big(}\tag{By \eqref{newmu}}\\
		&= u^\prime \cdot \phi_{a_k}(v^\prime) \cdot \big(\phi_{a_k}\circ \phi_{b_k}\big)(w^\prime) 
			\cdot a_kb_kc_k \phantom{\Big(}\tag{By Proposition \ref{hyp3}}\\
		&= u^\prime \cdot \phi_{a_k}\Big(v^\prime \cdot \phi_{b_k}(w^\prime)\Big)
			\cdot a_kb_kc_k \tag{By Proposition \ref{hyp1}}\\
		&= (u^\prime \cdot a_k) \cdot \Big(v^\prime \cdot \phi_{b_k}(w^\prime) \cdot b_kc_k\Big) \tag{By \eqref{newmu}}\\
		&= (u^\prime \cdot a_k) \cdot \Big(\big(v^\prime\cdot b_k\big) \cdot 
			\big(w^\prime \cdot c_k\big)\Big) \tag{By \eqref{newmu}}\\
		&= u \cdot (v \cdot w)\phantom{\Big(}
	\end{align*}
}

	\section{Conditions on Total Systems for Associativity of $r$-fold SDPs}

	\label{BracketAxiomsPhis}

In this section, we translate the elementary associativity conditions
${\bf A}[k,j,i]$ into conditions on the actions and brackets
	$$\phi_k^j : H_k \longrightarrow \mathrm{Maps}(H_j,H_j)$$
	$$[\cdot,\cdot]_{kj}^l : H_k \times H_j \longrightarrow H_l$$
constituting the total system $S^{(r)}$.

\bigskip

Let us begin by working out an explicit example for illustration. Take
the associativity condition ${\bf A}[4,3,2]$,  which 
is just the following condition for all $a_4 \in H_4$, $b_3 \in H_3$, $c_2 \in H_2$.
	$$a_4\cdot (b_3\cdot c_2) = (a_4\cdot b_3)\cdot c_2$$
	Working to evaluate the $\mu$-multiplications occurring in the
	left hand side (and making occasional use of part 1 of Proposition \ref{posspit} 
	to omit parentheses), one obtains
	\begin{align*}
		a_4\cdot (b_3\cdot c_2) &= \phi_{a_4}(b_3\cdot c_2)\cdot a_4 \\
		&= \phi_{a_4}\Big([b_3,c_2]_{32}^1\cdot{}^{\phi_3^2(b_3)}c_2\cdot b_3\Big)\cdot a_4
			\tag{By Def.\ref{BigMuDef}}\\
		&= \bigg(\phi_{a_4}\Big([b_3,c_2]_{32}^1\Big)\cdot \phi_{a_4}\Big({}^{\phi_3^2(b_3)}c_2\Big)\bigg) \cdot \phi_{a_4}(b_3) \cdot a_4
			\tag{By Def.\ref{extendedf}}
	\end{align*}
	\begin{align*}
		= \Big({}^{\phi_4^1(a_4)}[b_3,c_2]_{32}^1\cdot [a_4,{}^{\phi_3^2(b_3)}c_2]_{42}^1\cdot{}^{\phi_4^2(a_4)}\big({}^{\phi_3^2(b_3)}c_2\big)\Big)\cdot \\
		\Big([a_4,b_3]_{43}^1\cdot[a_4,b_3]_{43}^2\cdot{}^{\phi_4^3(a_4)}b_3\Big) \cdot a_4
			\tag{By Def.\ref{conjcommop}}
	\end{align*}
	\begin{align*}  = {}^{\phi_4^1(a_4)}[b_3,c_2]_{32}^1\cdot [a_4,{}^{\phi_3^2(b_3)}c_2]_{42}^1\cdot {}^{{}^{\phi_4^2(a_4)}\big({}^{\phi_3^2(b_3)}c_2\big)}[a_4,b_3]_{43}^1 \\
		\cdot {}^{\phi_4^2(a_4)}\big({}^{\phi_3^2(b_3)}c_2\big)\cdot[a_4,b_3]_{43}^2 \\
		\cdot {}^{\phi_4^3(a_4)}b_3  \\
		\cdot a_4\tag{By Def.\ref{BigMuDef}}
	\end{align*}
where in the last expression the 1st, 2nd, 3rd and 4th
components of the product have been written on separate lines. Working similarly with the
right hand side, one obtains the following, in which the components have
also been written on separate lines.
	\begin{align*}[a_4,b_3]^1_{43}\cdot{}^{\phi_2^1([a_4,b_3]^2_{43})}\Big({}^{{\phi_3^1({}^{\phi_4^3(a_4)}b_3)}}[a_4,c_2]^1_{42}\Big)\cdot{}^{\phi_2^1([a_4,b_3]^2_{43})}[{}^{\phi_4^3(a_4)}b_3,{}^{\phi_4^3(a_4)}c_2]^1_{32} \\
		\cdot[a_4,b_3]^2_{43}\cdot {}^{\phi_3^2({}^{\phi_4^3(a_4)}b_3)}\big({}^{\phi_4^2(a_4)}c_2\big) \\
		\cdot {}^{\phi_4^3(a_4)}b_3 \\
		\cdot a_4
	\end{align*}
By equating the 1st components, a condition on the total system
$S^{(r)}$ is obtained that will be denoted ${\bf A}[4,3,2;1]$. Similarly the
condition ${\bf A}[4,3,2;2]$ is obtained by equating the 2nd components.
The conditions ${\bf A}[4,3,2;3]$ and ${\bf A}[4,3,2;4]$ are vacuous since
the 3rd and 4th components are identical. Evidently the two conditions
${\bf A}[4,3,2;1]$ and ${\bf A}[4,3,2;2]$ taken together are equivalent
to the condition ${\bf A}[4,3,2]$.



\defn{
	Let the symbol ${\bf A}[k,j,i;l]$ denote the algebraic condition
	obtained by using the formulas for $\mu$
	to separately evaluate each side of ${\bf A}[k,j,i]$, that is,
	the equation
		\begin{equation*}
			a_k\cdot(b_j\cdot c_i) = (a_k\cdot b_j)\cdot c_i
		\end{equation*}
	and then equating the $l$th components of each side.
}

\begin{prp}\label{vacuouslgi}
	Fix $k > j > i$.  
	For $l > i$, the $l$th components of the left- and right-hand sides of ${\bf A}[k,j,i]$
	are identical.  Therefore the conditions ${\bf A}[k,j,i;l]$ with $l>i$ are vacuous.
\end{prp}
\prf{Proof.}{
	Evaluate the left-hand side of ${\bf A}[k,j,i]$ as follows.
		\begin{align*}
			a_k\cdot (b_j\cdot c_i) &= a_k\cdot \big(\phi_{b_j}(c_i)\cdot b_j\big) 
				\tag{By \eqref{tiny}}\\
			&= \phi_{a_k}\big(\phi_{b_j}(c_i)\cdot b_j\big)\cdot a_k
				\tag{By \eqref{tiny}} \\
			&= \Big(\big(\phi_{a_k}\circ \phi_{b_j}\big)(c_i)\cdot \phi_{a_k}(b_j)\Big)\cdot a_k
				\tag{By \eqref{spit}}
		\end{align*}
	Now writing $u$ for the first $i$ components of $\phi_{a_k}(b_j)$ and $b_{i+1}^\p, \ldots, b_j^\p$
	for its	remaining components, the above becomes
		$$\Big(\big(\phi_{a_k}\circ \phi_{b_j}\big)(c_i)\cdot u\Big)\cdot b_{i+1}^\p 
			\cdot\ldots\cdot b_j^\p\cdot a_k.$$

	Now evaluate the right-hand side of ${\bf A}[k,j,i]$, using the same notation for 
	the components of $\phi_{a_k}(b_j)$.
		\begin{align*}
			(a_k\cdot b_j)\cdot c_i &= \big(\phi_{a_k}(b_j)\cdot a_k\big)\cdot c_i \tag{By \eqref{tiny}}\\
			&= \phi_{a_k}(b_j)\cdot \phi_{a_k}(c_i) \cdot a_k \tag{By \eqref{newmu}} \\
			&= (u\cdot b_{i+1}^\p \cdot\ldots\cdot b_j^\p)\cdot \phi_{a_k}(c_i) \cdot a_k \\
			&= u\cdot \big(\phi_{b_{i+1}^\p}\circ\ldots\circ \phi_{b_j^\p}\circ \phi_{a_k}\big)(c_i) \cdot b_{i+1}^\p \cdot\ldots\cdot b_j^\p\cdot a_k \tag{By \eqref{newmu}}
		\end{align*}

	Comparing the two results, evidently their $l$th components coincide for $l>i$.
}

\bigskip

The following notational
simplifications help compensate for the complexity of formulas such as the two
sides of ${\bf A}[4,3,2;1]$ above. Since the domain of a bracket is identified by the
subscripts of the arguments, the subscript on the
brackets may be left off, that is
	$$[a_k,b_j]^i := [a_k,b_j]^i_{kj}.$$
Similarly, the action $\phi_k^j$ is identified both by its argument and
by what the resulting endomorphism is acting on, so it will be abbreviated by writing
	$${}^{a_k}b_j := {}^{\phi_k^j(a_k)}b_j.$$
When composing actions, even coming from
different groups $H_k$ and $H_j$, a dot is used to denote the
composition. Thus
	$${}^{a_k\cdot b_j}c_i := {}^{a_k}\big({}^{b_j}c_i\big)
		={}^{\phi_k^i(a_k)}\big({}^{\phi_j^i(b_j)}c_i\big).$$
(Elsewhere the dot in $a_k\cdot b_j$ denotes the multiplication $\mu$,
whereas in left superscripts it denotes composition of actions.)
Finally, assume that the subscripts $k,j,i$ are always used on the
symbols $a,b,c$ respectively, so the subscripts may be suppressed 
altogether, that is
	$$a := a_k, ~~~b := b_j, ~~~ c := c_i.$$

Using these conventions with ${\bf A}[4,3,2;2]$, for example, its form above 
takes on a leaner appearance as follows.
	\begin{equation*}
		{}^{a\cdot b}c\cdot[a,b]^2 = [a,b]^2\cdot {}^{{}^{a}b\cdot a}c
	\end{equation*}
All terms of such expressions have unambiguous meanings via the above
conventions, provided the indices in the symbol ${\bf A}[k,j,i;l]$ are
also specified.

\bigskip

Working out ${\bf A}[4,2,1;1]$ using these conventions, one discovers that
it has the same form as ${\bf A}[4,3,2;2]$, the only difference being in
the upper indices.
	$${}^{\phi_4^1(a_4)}\big({}^{\phi_2^1(b_2)}c_1\big)\cdot[a_4,b_2]_{42}^1 = 
		[a_4,b_2]^1_{42}\cdot {}^{\phi_2^1({}^{\phi_4^2(a_4)}b_2)}\big({}^{\phi_4^1(a_4)}c_1\big)$$
	$${}^{a\cdot b}c\cdot[a,b]^1 = [a,b]^1\cdot {}^{{}^{a}b\cdot a}c$$
When two conditions ${\bf A}[k,j,i;l]$, ${\bf A}[k^\p,j^\p,i^\p;l^\p]$
bear this relation to one another, that is, they have exactly the same
form except for changing upper indices, call them {\it formally
similar} or just {\it similar} and denote this situation as follows.
	$$\xymatrix{{\bf A}[k,j,i;l] \ar@{ ~ }[r] & {\bf A}[k^\p,j^\p,i^\p;l^\p]}$$
The next proposition gives general conditions under which formal
similarities appear.

\begin{prp}  For any $k > j > i \ge l$ the following formal similarities hold.
	$$\xymatrix{{\bf A}[k+1,j,i;l]\ar@{ ~ }[r]^-1&{\bf A}[k,j,i;l]\ar@{ ~ }[r]^-2 & {\bf A}[k+1,j+1,i+1;l+1]}$$
	$$\xymatrix{{\bf A}[k+1,k+1,i;l]\ar@{ ~ }[r]^-3&{\bf A}[k,k,i;l]\ar@{ ~ }[r]^-4 & {\bf A}[k+1,k+1,i+1;l+1]}$$
\end{prp}
\prf{Proof.}{
	The proofs of the similarities labelled 3 and 4 are essentially the same
	as the proofs of the similarities labelled 1 and 2, so we explain only
	these.

	For the first similarity relation, note that in complete generality the SDP
		$$G = H_1\rtimes\ldots\rtimes H_r$$
	contains sub-SDPs of the form
		$$G_{j,k} := H_1\rtimes\ldots\rtimes H_j \rtimes H_k$$
	whose total system $S_{j,k}^{(j+1)}$ is obtained by forgetting all data
	involving groups $H_i$ other than the ones shown. This produces a
	well-defined total system because,
	according to Definition
	\ref{foperator},
	pushing $h_k$ across $h_j$ leaves $h_k$ untouched and causes new terms
	only of index less than or equal to $j$ to be introduced. Now the
	components of the multiplication $\mu_{j,k}$ of $G_{j,k}$ are identical
	to the corresponding components of $\mu$, so the condition ${\bf A}[j+1,j,i;l]$
	for $G_{j,k}$ (with its last component $H_k$ reindexed as its $(j+1)$st
	component) is precisely the condition ${\bf A}[k,j,i;l]$ for $G$. Moreover,
	since any $(j+1)$-SDP arises as $G_{j,k}$ for some $G$ (for
	instance, fill out the other components of $G$ with trivial groups), the
	form of the condition ${\bf A}[j+1,j,i;l]$ is the same as the one
	occurring in the general case. Therefore, for fixed $j,i,l$, the
	form of ${\bf A}[k,j,i;l]$ is independent of $k$ for $k>j$.	

	Now for the second similarity relation.  Inspecting once again the 
	formula for the multiplication $\mu$ on an external SDP
		$$G := H_1\rtimes\ldots\rtimes H_r$$
	defined via a total system $S^{(r)}$, one sees that the components of
	$\mu$ in positions $2,\ldots r$ do not involve the group $H_1$. Hence one
	may form a new total system $S^{(r-1)}/H_1$, obtained from $S^{(r)}$ by
	forgetting all information involving $H_1$, and use it to define a
	quotient SDP
		$$G/H_1 = H_2 \rtimes\ldots\rtimes H_r$$
	with multiplication denoted $\mu/H_1$.  Note the projection map 
		$$G\lra G/H_1$$
		$$H_1\times\ldots\times H_r \lra H_2\times\ldots\times H_r$$
	is then a homomorphism with respect to $\mu$ and $\mu/H_1$. Since the
	formulas for the components $1,\ldots, r-1$ of $\mu/H_1$ are identical
	to those for components $2, \ldots, r$ respectively of $\mu$, the
	associativity condition ${\bf A}[k,j,i;l]$ for $\mu/H_1$ is just the
	associativity condition ${\bf A}[k+1,j+1,i+1;l+1]$ for $\mu$. Moreover, since
	any $(j-1)$-SDP arises as $G/H_1$ for some $G$ (for instance
	$\{1\}\rtimes H_2 \rtimes\ldots\rtimes H_r$), the form of the
	condition ${\bf A}[k,j,i;l]$ for $G/H_1$ is the same that occurring in the
	general case. Consequently, the condition ${\bf A}[k+1,j+1,i+1;l+1]$ for $G$
	has the same form as the condition ${\bf A}[k,j,i;l]$ occurring in the
	general case. 
}

We now use these similarity relations to select a maximal collection of
formally dissimilar conditions from among all associativity conditions
${\bf A}[k,j,i;l]$. Using similarities 2 and 4, any condition ${\bf
A}[k,j,i;l]$ is similar to one with $l=1$. Using similarity 1, any
condition ${\bf A}[k,j,i;1]$ with $k>j$ is similar to ${\bf
A}[j+1,j,i;1]$, and using similarity 3, any condition ${\bf A}[k,k,i;1]$
is similar to ${\bf A}[i+1,i+1,i;1]$. It follows that in each similarity
class there is a condition of one of the forms
	$${\bf A}[k, k-1, i;1] \text{~~or~~} {\bf A}[k,k,k-1;1]$$
where the indices have been reparametrized so that $k$ is the largest
index. Thus for each $k \ge 2$ one obtains $k-1$ different forms ${\bf
A}[k, k-1,i;1]$ corresponding to the values $1 \le i \le k-1$ and one
more form ${\bf A}[k,k,k-1;1]$ for a total of $k$ different forms of
conditions with largest index $k$. 

\bigskip

Here is a table of these forms for $k = 2,3,4,5$. The reader will be
able to translate the conditions for $k=2$ into the familiar ones for
the usual external semidirect product construction (external 2-SDP).

\begin{equation*} \begin{split}
	{\bf A}[2,1,1;1]&: \mathrm{Image}(\phi_2^1) \subseteq \mathrm{End}(H_1)\\
	{\bf A}[2,2,1;1]&: \phi_2^1~\text{is a homomorphism}\\[2pt]
	\hline\\[-16pt]
	{\bf A}[3,2,1;1]&: [a,b]^1\cdot {}^{{}^{a}b\cdot a}c = {}^{a\cdot b}c\cdot[a,b]^1\\
	{\bf A}[3,2,2;1]&: [a, bc]^1 = [a,b]^1\cdot{}^{{}^{a}b}[a,c]^1 \\
	{\bf A}[3,3,2;1]&: [ab,c]^1 = {}^{a}[b,c]^1\cdot[a,{}^{b}c]^1\\[2pt]
	\hline\\[-16pt]
	{\bf A}[4,3,1;1]&: [a,b]^1\cdot{}^{[a,b]^2\cdot{}^{a}b\cdot a}c = 
			{}^{a\cdot b}c\cdot[a,b]^1\\
	{\bf A}[4,3,2;1]&: [a,b]^1\cdot{}^{[a,b]^2\cdot{}^{a}b}[a,c]^1\cdot{}^{[a,b]^2}[{}^{a}b,{}^{a}c]^1 =
		{}^{a}[b,c]^1 \cdot [a,{}^{b}c]^1 \cdot {}^{{}^{{}^{a\cdot b}}c}[a,b]^1	\\
	{\bf A}[4,3,3;1]&: [a,bc]^1 = [a,b]^1\cdot{}^{[a,b]^2\cdot{}^{a}b}[a,c]^1\cdot
			{}^{[a,b]^2}[{}^{a}b,[a,c]^2]^1\\
	{\bf A}[4,4,3;1]&: [ab,c]^1 = {}^{a}[b,c]^1\cdot[a,[b,c]^2]^1\cdot{}^{{}^{a}[b,c]^2}[a,{}^{b}c]^1\\[2pt]
	\hline\\[-16pt]
	{\bf A}[5,4,1;1]&: [a,b]^1\cdot{}^{[a,b]^2\cdot[a,b]^3\cdot{}^{a}b\cdot a}c = {}^{a\cdot b}c\cdot[a,b]^1\\
	{\bf A}[5,4,2;1]&: [a,b]^1\cdot{}^{[a,b]^2\cdot[a,b]^3\cdot{}^{a}b}[a,c]^1\cdot{}^{[a,b]^3}[{}^{a}b,{}^{a}c]^1\cdot[[a,b]^3,{}^{{}^{a}b\cdot a}c]^1 = \\
		& \hspace{20pt} {}^{a}[b,c]^1\cdot[a,{}^{b}c]^1\cdot{}^{{}^{a\cdot b}c}[a,b]^1 \\
	{\bf A}[5,4,3;1]&: [a,b]^1\cdot{}^{[a,b]^2\cdot[a,b]^3\cdot{}^{a}b}[a,c]^1\cdot{}^{[a,b]^2\cdot[a,b]^3}[{}^{a}b,[a,c]^2]^1\cdot{}^{[a,b]^3}[[a,b]^3,{}^{{}^{a}b}[a,c]^2]^1\cdot \\
		& \hspace{20pt} {}^{[a,b]^2\cdot{}^{[a,b]^3\cdot{}^{a}b}[a,c]^2\cdot[a,b]^3}[{}^{a}b,{}^{a}c]^1\cdot{}^{[a,b]^2\cdot{}^{[a,b]^3\cdot{}^{a}b}[a,c]^2}[[a,b]^3,[{}^{a}b,{}^{a}c]^2]^1 = \\
		& \hspace{30pt} {}^{a}[b,c]^1\cdot[a,[b,c]^2]^1\cdot{}^{{}^{a}[b,c]^2}[a,{}^{b}c]^1\cdot 
			{}^{{}^{a}[b,c]^2\cdot[a,{}^{b}c]^2\cdot a\cdot{}^{b}c}[a,b]^1\cdot \\ 
		& \hspace{40pt}{}^{{}^{a}[b,c]^2\cdot[a,{}^{b}c]^2}[{}^{a\cdot b}c,[a,b]^2]^1 \\
	{\bf A}[5,4,4;1]&: [a,bc]^1 = [a,b]^1\cdot{}^{[a,b]^2\cdot[a,b]^3\cdot{}^{a}b}[a,c]^1\cdot[{}^{a}b,[a,c]^2]^1
				\cdot{}^{[a,b]^2}[[a,b]^3,{}^{{}^{a}b}[a,c]^2]^1 \\
		& \hspace{20pt} \cdot {}^{[a,b]^2\cdot{}^{[a,b]^3\cdot{}^{a}b}[a,c]^2\cdot[a,b]^3}[{}^{a}b,[a,c]^3]^1\cdot{}^{[a,b]^2\cdot{}^{[a,b]^3\cdot{}^{a}b}[a,c]^2}[[a,b]^3,[{}^{a}b,[a,c]^3]^2]^1\\
	{\bf A}[5,5,4;1]&: [ab,c]^1 = {}^{a}[b,c]^1\cdot[a,[b,c]^2]^1\cdot{}^{{}^{a}[b,c]^2}[a,[b,c]^3]^1 \\ 
		& \hspace{20pt} \cdot
			{}^{{}^{a}[b,c]^2\cdot[a,[b,c]^3]^2\cdot{}^{a}[b,c]^3}[a,{}^{b}c]^1\cdot{}^{{}^{a}[b,c]^2\cdot[a,[b,c]^3]^2}[{}^{a}[b,c]^3,[a,{}^{b}c]^2]^1 \\[2pt]
	\hline\\[-16pt]
\end{split}\end{equation*}

	\section{Conditions for Homomorphisms From an $r$-fold SDP}		

	\label{HomFromSDPPhis}

In this final section, we prove a result similar in spirit to that of
the section \ref{GenAssocRelsNewPhis}, to the effect that, if $G^\p$ is
associative (i.e., a semigroup), then in order to check that a map 
	$$f: H_1\rtimes\ldots\rtimes H_r \lra G^\p$$
is homomorphism, it suffices to check all relations of the form
	$$f(a_k\cdot b_j) = f(a_k)\cdot f(b_j)$$
where $a_k\in H_k$, $b_j\in H_j$ and $k \ge j$.
A result having the same purpose but relating to
homomorphisms {\it into} an SDP can be found in \cite{CarrascoCegarra}. 

\bigskip

Let $G = H_1\rtimes\ldots\rtimes H_r$, $\mu$, $S^{(r)}$, $R_k$ and the
operators $\phi_{a_k}$ be as in previous sections. Also let $G^\p$ with an
associative multiplication be given, along with maps
	$$f_i: H_i \lra G^\p$$
for each $i$.  These assemble into a (not necessarily homomorphic) map $f$ via the formula
	$$f: G \lra G^\p$$
	\begin{equation}\label{phicompdef}
		f(h_1\cdot\ldots\cdot h_r) := f(h_1)\cdot\ldots\cdot f(h_r)
	\end{equation}
in which
	$$f(h_i) := f_i(h_i)$$
for each $i$.  Also the dots $\cdot$
on the left-hand side denote the external SDP multiplication $\mu$,
while the dots on the right-hand side denote the multiplication in $G^\p$.

\defn{
	Let $U, V$ be subsets of $G$.  Let the symbol
		$${\bf H}[f; U,V]$$
	stand for the family of all homomorphism relations of the form
		$$f(u\cdot v) = f(u)\cdot f(v)$$
	for $u\in U, v\in V$.  
\newline\indent
	In the special case in which one of $U,V$ is $R_k$ or $H_k$, 
	the symbol $(k)$ or $k$ respectively will be substituted for it in the above notation.
	Here is	an example.
		$${\bf H}[f; k,(j)] := {\bf H}[f; H_k, R_j]$$
	We also use the following abbreviation.	
		$${\bf H}[[f; k]] := {\bf H}[f; R_k,R_k]$$
}

\bigskip

The main goal of this section is to prove the following result.

\begin{thm} \label{homsdpthm}
	Let $k$ with $r \ge k > 1$ be fixed. In order for the restriction of
	$f$ to $R_k$ to be a homomorphism, it suffices to require the
	following homomorphism conditions.
	\begin{enumerate}
		\item{${\bf H}[[f; k-1]]$, i.e., $f$ is a homomorphism on $R_{k-1}$.}
		\item{${\bf H}[f; k,j]$ for all $j$ with $k \ge j \ge 1$.}
	\end{enumerate}
\end{thm}

From this a simple induction argument yields the following.

\begin{thm}\label{hommaincor}
	In order for $f$ to be a homomorphism on $G$, it suffices to assume only that 
	the homomorphism conditions ${\bf H}[f; k,j]$ hold for all $j,k$ with
	$r \ge k \ge j \ge 1$.
\end{thm}


\bigskip

In the proof, use is again made of the identities \eqref{tiny} and
\eqref{spit} from Remark \ref{synopsis}. We add one more formula to the
list, which follows immediately from the formula for $f$ given above.
	\begin{align} 
		\label{homspit}
		f(u \cdot a_k) &= f(u)\cdot f(a_k) & \forall& a_k\in H_k, u\in R_{k-1}
	\end{align}

\begin{prp} \label{homhyp1}
	Let $k$ with $r \ge k > 1$ be fixed.  Assume ${\bf H}[f; [k-1]]$, and also
	assume ${\bf H}[f; k,j]$ for all $j$ such that $k > j \ge 1$.
	Then the condition 
		$${\bf H}[f; k,(k-1)]$$
	holds, or equivalently the formula
	\begin{equation*} 
		f(a_k\cdot v) = f(a_k)\cdot f(v)
	\end{equation*}
	is valid for any $a_k\in H_k$ and $v\in R_{k-1}$.
\end{prp}
\prf{Proof.}{
	We show that, for fixed $j$ with $k > j \ge 1$, under the assumptions 
		$${\bf H}[[f; k-1]] \textrm{~~and~~} {\bf H}[f; k,j]$$
	the following implication holds.
		$${\bf H}[f; k,(j-1)] \implies {\bf H}[f; k,(j)]$$
	The assertion of the proposition then follows by induction.

	Take $a_k \in H_k$ and $v\in R_j$ and write
		$$v = v^\prime \cdot b_j$$
	where $v^\prime \in R_{j-1}$ and $b_j \in H_j$ are determined by $v$.  Then compute as follows.
	\begin{align*}
		f(a_k\cdot v) &= f(a_k\cdot v^\p\cdot b_j)\phantom{\Big(}\\
		&= f\Big(\phi_{a_k}(v^\p\cdot b_j) \cdot a_k\Big) \tag{By \eqref{tiny}} \\
		&= f\Big(\phi_{a_k}(v^\p)\cdot \phi_{a_k}(b_j) \cdot a_k\Big) \tag{By \eqref{spit}} \\
		&= f\Big(\phi_{a_k}(v^\p)\cdot \phi_{a_k}(b_j)\Big) \cdot f(a_k) \tag{By \eqref{homspit}} \\
		&= f\Big(\phi_{a_k}(v^\p)\Big)\cdot f\Big(\phi_{a_k}(b_j)\Big) \cdot f(a_k) \tag{By ${\bf H}[[f; k-1]]$} \\
		&= f\Big(\phi_{a_k}(v^\p)\Big)\cdot f\Big(\phi_{a_k}(b_j) \cdot a_k\Big) \tag{By \eqref{homspit}} \\
		&= f\Big(\phi_{a_k}(v^\p)\Big)\cdot f(a_k\cdot b_j) \tag{By \eqref{tiny}} \\
		&= f\Big(\phi_{a_k}(v^\p)\Big)\cdot f(a_k)\cdot f(b_j) \tag{By ${\bf H}[f; k,j]$} \\
		&= f\Big(\phi_{a_k}(v^\p)\cdot a_k\Big)\cdot f(b_j) \tag{By \eqref{homspit}} \\
		&= f(a_k\cdot v^\p)\cdot f(b_j) \phantom{\Big(}\tag{By \eqref{tiny}} \\
		&= f(a_k)\cdot f(v^\p)\cdot f(b_j) \phantom{\Big(}\tag{By ${\bf H}[f; k,(j-1)]$} \\
		&= f(a_k)\cdot f(v^\p\cdot b_j) \phantom{\Big(}\tag{By \eqref{homspit}} \\
		&= f(a_k)\cdot f(v) \phantom{\Big(}
	\end{align*}
}

\prf{Proof of Theorem \ref{homsdpthm}.}{
	First observe that the hypotheses of the theorem include those of Proposition
	\ref{homhyp1} and so it holds, that is, one has \mbox{${\bf H}[f; k,(k-1)]$}. Now take $u,v \in R_k$
	and write
	\begin{align*}
		u &= u^\p\cdot a_k \\
		v &= v^\p\cdot b_k
	\end{align*}
	where $u^\prime, v^\prime \in R_{k-1}$ and $a_k, b_k \in H_k$ are
	determined by $u,v$ respectively.  Then compute as follows.  
	\begin{align*}
		f(u\cdot v) &= f\Big(u^\p \cdot \phi_{a_k}(v^\p) \cdot a_kb_k\Big) \tag{By \eqref{newmu}} \\
		&= f\Big(u^\p \cdot \phi_{a_k}(v^\p)\Big) \cdot f(a_kb_k) \tag{By \eqref{homspit}} \\
		&= f(u^\p) \cdot f\Big(\phi_{a_k}(v^\p)\Big) \cdot f(a_k)\cdot f(b_k) \tag{By ${\bf H}[[f; k-1]]$ and ${\bf H}[f; k,k]$} \\
		&= f(u^\p) \cdot f\Big(\phi_{a_k}(v^\p)\cdot a_k\Big)\cdot f(b_k) \tag{By \eqref{homspit}} \\
		&= f(u^\p) \cdot f(a_k\cdot v^\p)\cdot f(b_k) \phantom{\Big(}\tag{By \eqref{tiny}} \\
		&= f(u^\p) \cdot f(a_k)\cdot f(v^\p)\cdot f(b_k)\phantom{\Big(} \tag{By Prop. \ref{homhyp1}} \\
		&= f(u^\p\cdot a_k)\cdot f(v^\p\cdot b_k) \phantom{\Big(}\tag{By \eqref{homspit}} \\
		&= f(u)\cdot f(v)\phantom{\Big(}
	\end{align*}
}

\bigskip

%
%


The following corollary is useful for computing homomorphism
conditions explicitly in practice.

\begin{cor}\label{phiprescomms}
	If $G^\p$ and the SDP $H_1 \rtimes \ldots \rtimes H_r$ are groups and
	the components $f_i$ of $f$ are homomorphisms, then the map $f$
	as given by \eqref{phicompdef} is a homomorphism if and only if it
	preserves all commutators of the form
		$$f\big([a_k, b_j]\big) = \big[f(a_k), f(b_j)\big]$$
	where $a_k\in H_k$, $b_j\in H_j$ and $k > j$.  If $H_j$ is closed under conjugation 
	by $H_k$ for a particular pair of indices $k>j$, then this condition simplifies to the following.
		$$f\big({}^{a_k}b_j\big) = {}^{f(a_k)}f(b_j)$$
\end{cor}
\prf{Proof.}{
	The following preliminary identity is used to obtain the first assertion.
		\begin{align*}
			f(a_k\cdot b_j) &= f\Bigg(\prod_{i=0}^{j-1}[a_k,b_j]_{kj}^i \cdot {}^{\phi_k^j(a_k)}b_j\cdot a_k\Bigg) \\
			&= f\Bigg(\prod_{i=0}^{j-1}[a_k,b_j]_{kj}^i \cdot [a_k,b_j]_{kj}^j \cdot b_j\cdot a_k\Bigg) \\
			&= f\Bigg(\prod_{i=0}^{j-1}[a_k,b_j]_{kj}^i\Bigg) \cdot f\big([a_k,b_j]_{kj}^j \cdot b_j\big) \cdot f(a_k) \tag{By \eqref{phicompdef}}\\
			&= f\Bigg(\prod_{i=0}^{j-1}[a_k,b_j]_{kj}^i\Bigg) \cdot f\big([a_k,b_j]_{kj}^j\big) \cdot f(b_j) \cdot f(a_k) \tag{$f_j$ is a hom.}\\
			&= f\big([a_k,b_j]\big)\cdot f(b_j)\cdot f(a_k) \phantom{\Bigg(}\tag{By \eqref{phicompdef}}
		\end{align*}
	Now if $f$ preserves
	commutators as in the statement of the proposition, then the right-hand side becomes 
	$f(a_k)\cdot f(b_j)$ and the identity becomes ${\bf H}[f; k,j]$.  By Theorem \ref{hommaincor}
	one concludes that $f$ is a homomorphism.
\newline\indent
	For the second assertion, note that the condition ${}^{H_k}H_j \subseteq H_j$
	implies that the product $\prod_{i=0}^{j-1}[a_k,b_j]_{kj}^i$ is trivial, so that
		$$[a_k,b_j] = {}^{\phi_j^k(a_k)}b_jb_j^{-1}$$
	and hence the action $\phi_j^k$ coincides with conjugation in $G$.
	Since $f$ is a homomorphism on $H_j$, the condition in question becomes
		\begin{align*}
			f\big([a_k, b_j]\big) &= \big[f(a_k), f(b_j)\big]\\
			f\big({}^{a_k}b_j\big)f\big(b_j^{-1}\big) &= {}^{f(a_k)}f(b_j)f(b_j)^{-1}
		\end{align*}
		$$\iff f\big({}^{a_k}b_j\big) = {}^{f(a_k)}f(b_j)$$
	as claimed.
}

\bibliographystyle{amsalpha}
\bibliography{Preprint_1_SDPs}

\end{document}